\documentclass{amsart}

\theoremstyle{plain}
\newtheorem{thm}{\indent\sc Theorem}[section]
\newtheorem{prop}[thm]{\indent\sc Proposition}

\newtheorem{cor}[thm]{\indent\sc Corollary}
\theoremstyle{definition}
\newtheorem{defn}[thm]{\indent\sc Definition}

\theoremstyle{remark}

\def\proofwp{{\sc Proof}}

\numberwithin{equation}{section}
\def\cA{\mathcal{A}}

\def\cC{\mathcal{C}}

\def\cF{\mathcal{F}}

\usepackage{amscd}
\usepackage{amsmath,amssymb}
\usepackage{amssymb}
\usepackage{amsmath}
\usepackage{amsthm}
\usepackage{latexsym}
\usepackage{amsfonts}
\usepackage{color}
\usepackage{graphics}
\usepackage[all,knot]{xy}
\usepackage[all]{xy}
\xyoption{arc}
\usepackage{enumerate}
\usepackage{amstext}

\newcommand{\RR}{\mathbb{R}}

\begin{document}

\title[Rham cohomology of locally trivial Lie groupoids]{On Rham cohomology of locally trivial Lie groupoids over triangulated manifolds}



\author{{\textsc{Jose R. Oliveira}}}
\address{Department of Mathematics, Minho University,
    Braga, Portugal}
\curraddr{}
\email{jmo@math.uminho.pt}
\thanks{The author is partially supported by MICINN, Grant MTM2014-56950-P}



\date{}

\dedicatory{}


\maketitle






\begin{center}

\vspace{3mm}

\textsc{Abstract}

\vspace{3mm}

\end{center}


Based in the isomorphism between Lie algebroid cohomology and piecewise smooth cohomology, it is proved that the Rham cohomology of a locally trivial Lie groupoid $G$ on a smooth manifold $M$ is isomorphic to the piecewise Rham cohomology of $G$, in which $G$ and $M$ are manifolds without boundary and $M$ is smoothly triangulated by a finite simplicial complex $K$ such that, for each simplex $\Delta$ of $K$, the inverse images of $\Delta$ by the source and target mappings of $G$ are transverses submanifolds in the ambient space $G$. As a consequence, it is shown that the piecewise de Rham cohomology of $G$ does not depend on the triangulation of the base.




\vspace{3mm}







\section{Introduction}

\vspace{3mm}

It is well known that the de Rham theorem states that the integration map over smooth singular chains on a smooth manifold induces an isomorphism between the de Rham cohomology and the singular cohomology of that manifold (see \cite{weil-and}).

When a manifold is smoothly triangulated by a simplicial complex, differential form can be defined on the corresponding simplicial space as being a family of differential forms defined on the simplices of the simplicial complex which agree on the intersection of two simplices and satisfy a regular property when restricted to each simplex. A differential can be defined by the corresponding exterior derivative on each simplex, yielding a cochain algebra. Sullivan in \cite{suli-inf} and Whitney in \cite{wity-git} show the cohomology of that cochain algebra is isomorphic to classic cohomologies of the space.

The Rham-Sullivan theorem (see Sullivan \cite{suli-inf}, theorem 7.1) implies that the restriction map from the cochain algebra of the smooth forms on the manifold to the cochain algebra of the piecewise smooth forms on the simplicial space induces an isomorphism in cohomology.

Mishchenko and Oliveira in \cite{mish-oli} have extended similar those constructions considered by Whitney and Sullivan to transitive Lie algebroids defined over triangulated manifolds in which the regularity property of differential forms is the smoothness. They proved that the restriction map induces an isomorphism in cohomology (for details, see \cite{mish-oli}).

Since Lie algebroids are the infinitesimal objects of Lie groupoids, the present work is an extension of the Mishchenko-Oliveira's work to Lie groupoids defined over triangulated manifolds. In order to clarify better the intrinsic nature of this extension, fix a locally trivial Lie groupoid $G$ on a closed smooth manifold $M$, with $G$ also without boundary, and a triangulation of $M$ such that the inverse images of each simplex $\Delta$ by the source and target mappings of $G$ are transverses submanifolds in $G$. Since $G$ is locally trivial, its Lie algebroid is transitive and the Lie algebroid of the restriction of $G$ to a simplex coincides with the Lie algebroid of $G$ restricted to that simplex. Likewise to the notion of piecewise smooth forms on Lie algebroids defined by Mishchenko-Oliveira in \cite{mish-oli}, an algebra of invariant piecewise smooth forms can be defined on $G$. As shown in \cite{kuki-pra} or \cite{Wein-Xu}, the de Rham algebra made of all invariant smooth forms on $G$ is isomorphic to the algebra of all smooth forms of its Lie algebroid. This isomorphism jointly with the isomorphism given in Mishchenko-Oliveira's theorem (theorem 5.1 of \cite{mish-oli}) will ensure that the algebra made of all invariant piecewise smooth forms on $G$ is quasi-isomorphic to the algebra of all piecewise smooth forms of its Lie algebroid. We will use this quasi-isomorphism and the corollary 5.5 of \cite{mish-oli} to show that the piecewise de Rham cohomology of $G$ does not depend on the triangulation of $M$.




\vspace{3mm}

\textbf{Acknowledgments}. I want to thank to Aleksandr Mishchenko, Jesus Alvarez, James Stasheff and Nicolae Teleman for their strong dynamism to discuss several topics concerning this work.


\vspace{6mm}

\begin{center}
\section{Rham cohomology of Lie groupoids}
\end{center}

\vspace{3mm}

We begin by reviewing briefly basic definitions and constructions concerning Lie groupoids and Lie algebroids. An extensive discussion on these issues can be found in the book \cite{makz-lga} by Mackenzie. The papers \cite{kuki-prap}, \cite{kuki-chw}, \cite{makz-acla}, \cite{mish-cla}, \cite{mish-latb}, \cite{mish-obs} and \cite{jean-prad} contain a detailed exposition and many examples. Various definitions and properties stated through the entire paper can be found, on level of cell spaces, in \cite{jesus-cal}, \cite{suli-inf}, \cite{suli-Tok} and \cite{wity-git}. Throughout this paper, we shall work on manifolds which are smooth, finite-dimensional and possibly with boundaries of different indices.

\vspace{3mm}

\textit{Lie groupoids}. Let $G$ and $M$ be two smooth manifolds. The manifold $G$ is called a Lie groupoid with base $M$ if the following is given: two surjective submersions $\alpha:G\longrightarrow M$ and $\beta:G\longrightarrow M$, called the source projection and the target projection respectively, a smooth mapping $1:M\longrightarrow G$ called the object inclusion mapping, a smooth multiplication in $G\ast G=\{(h,g)\in G\times G: \alpha(h)=\beta(g)\}$ and a mapping $G\longrightarrow G$ called inverse mapping and denoted by $g\longrightarrow g^{-1}$, satisfying the identities shown in the definitions 1.1.1 and 1.1.3 of Mackenzie's book \cite{makz-lga}.

The smooth mapping $$(\beta,\alpha):G\longrightarrow M\times M$$ $$\ \ \ \ \ \ g\longrightarrow (\beta(g),\alpha(g))$$ is called the anchor of the Lie groupoid $G$ with base $M$. The Lie groupoid $G$ is said to be \textit{locally trivial} if its anchor is a surjective submersion. For each $x, z\in M$, the manifolds $G_{x}=\alpha^{-1}(x)$ and $G^{z}=\beta^{-1}(z)$ are called the $\alpha$-fibre of G in $x$ and the $\beta$-fibre of G in $z$ respectively. Analogously, if $X$ and $Z$ are subsets of $M$, the sets $\alpha^{-1}(X)$ and $\beta^{-1}(Z)$ of $G$ are denoted by $G_{X}$ and $G^{Z}$ respectively.

We assume now that $G$ and $M$ are manifolds without boundary. If $X$ and $Z$ are submanifolds of $M$, possibly with boundaries of different indices, then the sets $G_{X}$ and $G^{Z}$ are submanifolds of $G$. Let $N$ be a submanifold of $M$, possibly with boundaries of different indices, such that the manifolds $G_{N}$ and $G^{N}$ are transverses submanifolds in the ambient space $G$. By taking the restrictions of the projections and of the object inclusion mapping, the manifold $G_{N}\cap G^{N}$ is a Lie groupoid with base $N$. The Lie groupoid $G_{N}\cap G^{N}$ will be called the Lie groupoid restriction of $G$ to $N$ and denoted by $G^{!!}_{N}$.

\vspace{3mm}

\textit{Lie algebroids}. Let $M$ be a smooth manifold, $TM$ the tangent bundle to $M$ and $\Gamma(TM)$ the Lie algebra of the vector fields on $M$. A Lie algebroid on $M$ is a vector bundle $\pi:\cA\longrightarrow M$ with base $M$ equipped with a vector bundle morphism $\rho:\cA\longrightarrow TM$, called anchor of $\cA$, and a structure of real Lie algebra on the vector space $\Gamma(\cA)$ of the sections of $\cA$ such that the map $\rho_{\Gamma}:\Gamma(\cA)\longrightarrow\Gamma(TM)$, induced by $\rho$, is a Lie algebra homomorphism and the action of the algebra $\cC^{\infty}(M)$ on $\Gamma(\cA)$ satisfies the natural condition: $$[\xi,f\eta]=f[\xi,\eta] + (\rho_{\Gamma}(\xi)\cdot f)\eta$$ for each $\xi$, $\eta$ $\in \Gamma(\cA)$ and $f\in \cC^{\infty}(M)$. The Lie algebroid $\cA$ is called transitive if the anchor $\gamma$ is surjective.

Let $\varphi:N\hookrightarrow M$ be a submanifold, possibly with boundaries of different indices and assume that $\cA$ is transitive. We recall that the Lie algebroid restriction of $\cA$ to the submanifold $N$, denoted by $\cA^{!!}_{N}$, is the Lie algebroid $\varphi^{!!}\cA$ constructed as inverse image of $\cA$ by the mapping $\varphi$ (see \cite{makz-lga}, \cite{makz-acla} and \cite{mish-oli} for more details).

\vspace{3mm}

\textit{From Lie groupoids to Lie algebroids}. We recall the construction of the Lie algebroid of a Lie groupoid. Let $M$ be a smooth manifold and $G$ a Lie groupoid on $M$ with source projection $\alpha:G\longrightarrow M$ and target projection $\beta:G\longrightarrow M$. Denote by $1:M\longrightarrow G$ the object inclusion mapping of $G$ and $G_{x}=\alpha^{-1}(x)$ the $\alpha$-fibre of G in $x$, for each $x\in M$. Denote by $\cA(G)$ the disjoint union $\bigsqcup _{x\in M}T_{1_{x}}G_{x}$ equipped with the structure of vector bundle on $M$. Consider the map $\rho:\cA(G)\longrightarrow TM$ defined by $\rho(a)=D\beta_{1_{x}}(a)$. Define a Lie bracket on $\Gamma(\cA(G))$ in the following way: for each $\xi$ and $\eta$ $\in \Gamma(\cA(G))$, the Lie bracket is defined by $$[\xi,\eta]=[\xi',\eta']_{G}$$ in which $\xi'$ and $\eta'$ denote the unique $\alpha$-right-invariant vector fields on $G$ such that $\xi'_{1_{x}}=\xi_{x}$ and $\eta'_{1_{x}}=\eta_{x}$, $\forall x\in M$ (see Kubarski \cite{kuki-pra}). Then, $(\cA(G),[\cdot,\cdot],\rho)$ is a Lie algebroid on $M$ and is called the Lie algebroid of the Lie groupoid $G$.

Lie groupoids and Lie algebroids enjoy some of the properties of Lie groups and Lie algebras. We notice that not every Lie algebroid is integrable to a Lie groupoid. The theorem 4.1 of the paper \cite{crac-loja} shows necessary and sufficient conditions so that a Lie algebroid is integrable to a Lie groupoid.

The proposition 3.5.18 of \cite{makz-lga} sates that, if the Lie groupoid $G$ is locally trivial, the Lie algebroid $\cA(G)$ is transitive.

Next proposition relates the restriction of a locally trivial Lie groupoid with the restriction of its Lie algebroid.

\vspace{3mm}

\begin{prop} Let $G$ be a locally trivial Lie groupoid on a smooth manifold $M$, in which $G$ and $M$ are manifolds without boundary, and $\cA(G)$ denote the Lie algebroid of $G$. Let $N$ be a submanifold of $M$ such that $G_{N}$ and $G^{N}$ are transverses submanifolds in the ambient space $G$. Then, the Lie algebroid of $G^{!!}_{N}$ is the Lie algebroid $\cA^{!!}_{N}$.
\end{prop}





\textit{Smooth forms on Lie groupoids}. Let $G$ be a Lie groupoid on a smooth manifold $M$ with source projection $\alpha:G\longrightarrow M$, target projection $\beta:G\longrightarrow M$ and object inclusion mapping $1:M\longrightarrow G$. Since $\alpha$ is a surjective submersion, it induces a foliation $\cF$ on $G$. Let $T \cF$ denote the tangent bundle of $\cF$ and $\RR_{M}=M\times \RR$ the trivial vector bundle on $M$ of fibre $\RR$.


\begin{defn} A smooth $\alpha$-form of degree $p$ on the Lie groupoid $G$ is a smooth section of the exterior vector bundle $\bigwedge ^{p}(T^{\ast}\cF; \RR_{M})$.
\end{defn}

Therefore, a smooth form on $G$ is a family $\omega=(\omega_{g})_{g\in G}$ such that, for each $g\in G$, one has $$\omega_{g}\in \bigwedge^{p}\big(\bigsqcup _{g\in G}T_{g}^{\ast}G_{\alpha(g)};\RR\big)$$ The set $\Omega_{\alpha}^{\ast}(G;M)$ of all smooth $\alpha$-forms on $G$ is a commutative graded algebra. The usual exterior derivative along the $\alpha$-fibres is defined by $$(d^{p}_{\alpha}\omega)(X_{1},X_{2},\cdot\cdot\cdot,X_{p+1})=\sum_{j=1}^{p+1}(-1)^{j+1}X_{j}\cdot(\omega
(X_{1},\cdot\cdot\cdot,\widehat{X_{j}},\cdot\cdot\cdot,X_{p+1}))\
+$$
$$+\sum_{i<k}(-1)^{i+k}\omega([X_{i},X_{k}],X_{1},\cdot\cdot\cdot,\widehat{X_{i}},\cdot\cdot\cdot,\widehat{X_{k}},\cdot\cdot\cdot,X_{p+1})$$ in which $\omega \in \Omega_{\alpha}^{p}(G;M)$ and $X_{1}$, $X_{2}$, $\cdot\cdot\cdot$, $X_{p+1}$ are smooth vector $\alpha$-fields on $G$. The complex $\Omega_{\alpha}^{\ast}(G;M)$ is a cochain algebra defined over $\RR$.

We are going now to consider invariant forms. For each $g\in G$, the right translation $R_{g}$ corresponding to $g$ is the mapping $R_{g}:G_{\beta(g)}\longrightarrow G_{\alpha(g)}$ defined by $R_{g}(h)=hg$. A smooth form $\omega\in \Omega_{\alpha}^{\ast}(G;M)$ is called right invariant or simply invariant if the equality
$$(R_{g})^{\ast}((\iota_{\alpha(g)})^{\ast}\omega)=(\iota_{\beta(g)})^{\ast}\omega$$ holds for each $g\in G$, in which $\iota_{x}:G_{x}\hookrightarrow G$ denotes the inclusion mapping for each $x\in M$ (see the third section of \cite{kuki-pra} in which the representation is the trivial representation on the trivial bundle $M\times \RR$). The set $\Omega^{\ast}_{\alpha,R}(G;M)$, consisting of all $\alpha$-forms on $G$ which are invariant under all groupoid right translations, is a subcomplex of $(\Omega_{\alpha}^{\ast}(G;M),d^{\ast}_{\alpha})$ and hence it is a cochain algebra.

\vspace{3mm}

\begin{defn} Keeping the same hypotheses and notations as above, the Rham cohomology of $G$ is the cohomology vector space of the cochain algebra $\Omega^{\ast}_{\alpha,R}(G;M)$. This cohomology vector space is denoted by $H^{\ast}_{\alpha,R}(G;M)$.
\end{defn}

\vspace{3mm}

\textit{Smooth forms on Lie algebroids}. We shall recall briefly the notion of smooth forms on Lie algebroids and its cohomology. Let $M$ be a smooth manifold and $\cA$ a Lie algebroid on $M$. A smooth form on $\cA$ is a section of $\bigwedge^{\ast}\big(\cA^{\ast};M\times \RR\big)$. The set $\Omega^{\ast}(\cA;M)$ of all smooth forms on $\cA$ is a commutative cochain algebra defined over $\RR$, in which the differential is given by $$d^{p}:\Omega^{p}(\cA;M)\longrightarrow \Omega^{p+1}(\cA;M)$$ $$d^{p}\omega (X_{1},X_{2},\cdot\cdot\cdot,X_{p+1})=\sum_{j=1}^{p+1}(-1)^{j+1}(\gamma\circ
X_{j})\cdot(\omega
(X_{1},\cdot\cdot\cdot,\widehat{X_{j}},\cdot\cdot\cdot,X_{p+1}))\
+$$
$$+\sum_{i<k}(-1)^{i+k}\omega([X_{i},X_{k}],X_{1},\cdot\cdot\cdot,\widehat{X_{i}},\cdot\cdot\cdot,\widehat{X_{k}},\cdot\cdot\cdot,X_{p+1})$$
for $\omega \in \Omega^{p}(\cA;M)$ and
$X_{1},X_{2},\cdot\cdot\cdot,X_{p+1}\in \Gamma (\cA)$.

The Lie algebroid cohomology of $\cA$ is the cohomology vector space of the cochain algebra $\Omega^{\ast}(\cA;M)$. This cohomology vector space is denoted by $H^{\ast}(\cA;M)$.



We notice that the proposition 7 of the paper \cite{kuki-pra}, combined with the subsequent paragraphs, states that the mapping
$$\Psi:\Omega^{\ast}_{\alpha,R}(G;M)\longrightarrow \Omega^{\ast}(\cA;M)$$ $$\Psi(\omega)_{x}=\omega_{1_{x}}$$ is an isomorphism of cochain algebras (cf. the paragraph preceding the theorem 1.2 of \cite{Wein-Xu}). In addition, we have $$H^{\ast}_{\alpha,R}(G;M)\cong H^{\ast}(\cA;M)$$ (cf. theorem 1.2 of \cite{Wein-Xu}).

\vspace{3mm}

\begin{center}
\section{Piecewise de Rham cohomology of Lie groupoids}
\end{center}

\vspace{3mm}

We begin by recalling Mishchenko-Oliveira's theorem (see \cite{mish-oli}, theorem 5.1). Let $M$ be a compact smooth manifold, smoothly triangulated by a simplicial complex $K$, and $\cA$ a transitive Lie algebroid on $K$. A piecewise smooth form on $\cA$ is a family $(\omega_{\Delta})_{\Delta \in K}$ of smooth forms, each form $\omega_{\Delta}$ defined on $\cA_{\Delta}^{!!}$, such that, if $\Delta$ and $\Delta'$ are simplices of $K$, with $\Delta'$ face of $\Delta$, $$(\omega_{\Delta})_{/\Delta'}=\omega_{\Delta'}$$ Let $\Omega^{\ast}(\cA;M)$ and $\Omega^{\ast}(\cA;K)$ denote respectively the cochain algebras of smooth forms and piecewise smooth forms on $\cA$. Consider the restriction mapping $$\Omega^{\ast}(\cA;M)\longrightarrow \Omega^{\ast}(\cA;K)$$ $$\omega \longrightarrow (\omega_{/\Delta})_{\Delta\in
K}$$ Mishchenko-Oliveira's theorem states that this mapping induces an isomorphism in cohomology.





We introduce now the notion of piecewise de Rham cohomology of Lie groupoids. Let $G$ be a locally trivial Lie groupoid with base $M$ such that $G$ and $M$ are manifolds without boundary and $M$ is smoothly triangulated by a simplicial complex $K$. Let $\alpha:G\longrightarrow M$, $\beta:G\longrightarrow M$ and $1:M\longrightarrow G$ denote respectively the source projection, the target projection and the object inclusion mapping. Suppose that, for each simplex $\Delta$ of $K$, the manifolds $G_{\Delta}$ and $G^{\Delta}$ are transverses submanifolds in the ambient space $G$. Analogous to piecewise smooth forms on Lie algebroids, we give now the notion of invariant piecewise form on $G$.

\vspace{3mm}

\begin{defn}
An invariant piecewise $\alpha$-form of degree $p$ on $G$ is a family $\omega=(\omega_{\Delta})_{\Delta \in K}$ such that, for each simplex $\Delta\in K$, $\omega_{\Delta}\in \Omega^{p}_{\alpha,R}(G^{!!}_{\Delta};\Delta)$ is an invariant $\alpha$-form of degree $p$ on $G^{!!}_{\Delta}$ and, if $\Delta$ and $\Delta'$ are two simplices of $K$ in which $\Delta'\prec \Delta$, the equality $(\omega_{\Delta})_{/\Delta'}=\omega_{\Delta'}$ holds.
\end{defn}

The $C^{\infty}(G)$-module of all invariant piecewise $\alpha$-forms of degree $p$ on $G$ is denoted by $\Omega^{p}_{\alpha,R}(G;K)$. A wedge product and an exterior derivative can be defined on the module $$\Omega^{\ast}_{\alpha,R}(G;K)=\bigoplus_{p\geq 0}\Omega^{p}_{\alpha,R}(G;K)$$ by the corresponding operations on each algebra $\Omega^{\ast}_{\alpha,R}(G^{!!}_{\Delta};{\Delta})$, giving to $\Omega^{\ast}_{\alpha,R}(G;K)$ a structure of cochain algebra defined over $\RR$.


\begin{defn} The piecewise de Rham cohomology of the Lie groupoid $G$, denoted by $H^{\ast}_{\alpha,R}(G;K)$, is the cohomology of the cochain algebra $\Omega^{\ast}_{\alpha,R}(G;K)$.
\end{defn}


Our aim is to relate the Rham cohomology $H^{\ast}_{\alpha,R}(G;M)$ of $G$ to the piecewise de Rham cohomology $H^{\ast}_{\alpha,R}(G;K)$ of $G$. For that, we have to consider a mapping $\phi$ from the complex $\Omega^{\ast}_{\alpha,R}(G;K)$ to the complex $\Omega^{\ast}(\cA(G);K)$. In order to obtain such mapping $\phi$, we recall that, for each simplex $\Delta$ of $K$, we have an isomorphism $$\psi_{\Delta}:\Omega^{p}_{\alpha,R}(G^{!!}_{\Delta};\Delta)\longrightarrow \Omega^{p}(\cA(G^{!!}_{\Delta});\Delta)$$ given by $(\psi_{\Delta}(\omega))_{x}=\omega_{1_{x}}$. Consider now an invariant piecewise $\alpha$-form $$\omega=(\omega_{\Delta})_{\Delta \in K} \in \Omega^{p}_{\alpha,R}(G;K)$$ For each simplex $\Delta\in K$, take the smooth form $\xi_{\Delta}=\psi_{\Delta}(\omega_{\Delta})\in \Omega^{p}(\cA(G^{!!}_{\Delta});\Delta)$. If $\Delta'$ is a simplex of $K$ such that $\Delta'$ is a face of $\Delta$, then $(\xi_{\Delta})_{/\Delta'}=\xi_{\Delta'}$ and so $\xi=(\xi_{\Delta})_{\Delta\in K}$ is a piecewise smooth form on $\cA(G)$. Keeping these hypotheses and notations, we state our next proposition.

\vspace{3mm}

\begin{prop} \label{teoquasikuki} The mapping $\Phi:\Omega^{\ast}_{\alpha,R}(G;K)\longrightarrow \Omega^{\ast}(\cA(G);K)$ defined by $$\Phi((\omega_{\Delta})_{\Delta \in K})=(\psi_{\Delta}(\omega_{\Delta}))_{\Delta\in K}$$ is well defined and is an isomorphism of cochain algebras.
\end{prop}

\vspace{3mm}

We can state now the main proposition of this paper. Denote by $r_{G}$ the restriction mapping $$r_{G}:\Omega^{\ast}_{\alpha,R}(G;M)\longrightarrow \Omega^{\ast}_{\alpha,R}(G;K)$$ $$r_{G}(\omega)=(\omega_{/\Delta})_{\Delta\in K}$$ and assume that $M$ is compact. Our proposition is the following.

\vspace{3mm}

\begin{prop} \label{mainteo2} The mapping $r_{G}$ induces an isomorphism in cohomology.
\end{prop}

\proofwp. The diagram
$$
\xymatrix{
\Omega^{\ast}_{\alpha,R}(G;M)\ar[d]_{r_{G}}\ar[r]^{iso} & \Omega^{\ast}(\cA;M)\ar[d]^{r_{\cA}} \\
 \Omega^{\ast}_{\alpha,R}(G;K)\ar[r]^{iso} & \Omega^{\ast}(\cA;K)
}
$$
is commutative, where $r_{\cA}$ is the restriction map given at the Mishchenko-Oliveira's theorem (see theorem 5.1 of \cite{mish-oli}). We apply the Mishchenko-Oliveira theorem 5.1 of \cite{mish-oli} in cohomology and the proof is done. {\small $\square$}

\vspace{3mm}

Our last proposition states that the piecewise de Rham cohomology of a locally trivial Lie groupoid on a compact triangulated manifold does not depend on the triangulation of the base, that is, for any simplicial subdivision of the simplicial complex, the piecewise de Rham cohomologies of the two triangulated manifolds are isomorphic. Precisely, this statement is our next proposition.

\vspace{3mm}

\begin{cor} Keeping the same hypotheses and notations as above, let $L$ be other simplicial complex which is a simplicial subdivision of $K$. Then, the piecewise de Rham cohomology of $G$ obtained by the triangulation corresponding to $K$ is isomorphic to the piecewise de Rham cohomology of $G$ obtained by the triangulation corresponding to $L$. In addition, this isomorphism is induced by the restriction map.
\end{cor}

\proofwp. Denote by $\phi:\Omega^{\ast}_{\alpha,R}(G;K)\longrightarrow \Omega^{\ast}_{\alpha,R}(G;L)$ the map given by restriction. The diagram
$$
\xymatrix{
& \Omega^{\ast}_{\alpha,R}(G;M) \ar[dl] \ar[dr] & & \\
\Omega^{\ast}_{\alpha,R}(G;K) \ar[rr]^{\phi} \ar[d] & & \Omega^{\ast}_{\alpha,R}(G;L)\ar[d] & & \\
\Omega^{p}(\cA;K) \ar[rr]^{\phi^{\cA}} & & \Omega^{p}(\cA;L)
}
$$
is commutative. By the corollary 5.5 of \cite{mish-oli}, the mapping $\phi^{\cA}$ induces an isomorphism in cohomology. By the propositions \ref{teoquasikuki} and \ref{mainteo2} above, the mappings non labeled induce isomorphisms in cohomology and so the mapping $\phi$ also induces an isomorphism in cohomology.


\end{document}